\documentclass[11pt,leqno,oneside,letterpaper]{amsart}
\usepackage[letterpaper,left=1.75truein, top=1truein]{geometry}
\usepackage{amsmath,amsfonts,amssymb,amscd,amsthm,amsbsy,epsf,graphics}

\usepackage{color}

\usepackage[colorlinks,linkcolor=blue,citecolor=blue]{hyperref}
\usepackage{epsfig}
\usepackage{graphicx}

\newtheorem{thm}{Theorem}
\newtheorem{cor}[thm]{Corollary}

\newtheorem{conj}[thm]{Conjecture}

\newtheorem{que}[thm]{Question}

\theoremstyle{definition}

\def\qed{{\hspace{2mm}{\small $\diamondsuit$}}}

\def\D{{\mathcal D}}

 \def\e{{\epsilon}}
 \def\l{{\lambda}}
 
 \def\m{{\mu}}

   \def\s{{\sigma}}
 
 \def\a{{\alpha}}
 \def\b{{\beta}}
 \def\p{{\partial}}
 \def\r{{\rho}}
 \def\ra{{\rightarrow}}

 \def\g{{\gamma}}
 \def\D{{\Delta}}

 \def\c{{\mathbb C}}
 \def\z{{\mathbb Z}}
 \def\2{{\mathbb Z_2}}

 \def\sl2{{SL(2,\mathbb C)}}
 
 \def\qed{{\hspace{2mm}{\small $\diamondsuit$}}}

 \def\pf{{\noindent{\bf Proof.\hspace{2mm}}}}

  \def\sl{{{\mbox{\footnotesize  $\mathfrak{L}$}}}}

\begin{document}

\title{A note on orderability and Dehn filling}

\author[Christopher Herald]{Christopher Herald}
\address{Department of Mathematics and Statistics,
 University of Nevada,
 Reno, NV 89557}
\email{herald@unr.edu}

\author{Xingru Zhang}
\address{Department of Mathematics, University at Buffalo, Buffalo, NY 14260}
\email{xinzhang@buffalo.edu}

\maketitle
\vspace{-.6cm}
\begin{center}
\today
\end{center}

\begin{abstract}
We improve upon a recent result of Culler and Dunfield  on orderability
of certain Dehn fillings by removing a difficult condition they
required.
\end{abstract}

If  $M$ is  the exterior of a knot $K$ in an integral homology $3$-sphere,  $\m$ and $\l$
will be the canonical meridian and longitude in
$\p M$, i.e.,  $\m$ bounds a meridian disk of $K$ and $\l$ is null-homologous
in $M$.  The set of slopes in $\p M$ will be identified  with
$\mathbb Q\cup \{1/0\}$ with respect to the chosen meridian and longitude in the usual way so that $\m$ is identified with $1/0$ and $\l$ with $0$,
and $M(r)$ will denote the Dehn filling
of $M$ with slope $r$.
Recall that a group   is called left-orderable if there is a total ordering on the group which is
invariant under left multiplication.
The purpose of this note is to update a recent result of Culler and Dunfield \cite{CD} to the following

\begin{thm}\label{new thm}Let $M$ be the exterior of a knot in an integral
homology $3$-sphere such that $M$ is irreducible.
If the Alexander polynomial $\D(t)$ of $M$ has a simple root on the unit
circle, then there exists a real number  $a>0$ such that, for every rational slope $r\in (-a,0)\cup(0,a)$, the Dehn filling $M(r)$ has left-orderable fundamental group.
\end{thm}

The above theorem was proved in \cite{CD} under
the additional condition that every closed essential surface in $M(0)$ is a fiber in a fibration of $M(0)$ over $S^1$ (\cite[Theorem 1.2]{CD}) or under the additional condition that each positive dimensional
component of the $PSL_2(\c)$ character variety of $M(0)$
consists entirely of characters of reducible representations \cite[Theorem 7.1]{CD}. The former  condition is  very restrictive and the latter
one is  hard to verify in general. So the updated theorem will be much more applicable.

We also remark that if $M(0)$ is prime (which is always true
if $M$ is the exterior of a   knot in $S^3$) then $\pi_1(M(0))$ is left-orderable since the first Betti number of $M(0)$ is positive (\cite[Theorem 1.3]{BRW}).
 In this case,  we may replace the intervals $(-a,0)\cup(0,a)$ in Theorem \ref{new thm} by the interval $(-a,a)$.

The motivation for studying if a Dehn filling has left-orderable fundamental group is its connection to the following now well known

\begin{conj}\label{conj}{\rm (\cite{BGW})} For a closed connected  orientable  prime $3$-manifold $W$, the following are
equivalent:\newline
(1) $W$ has left-orderable fundamental group.\newline
(2) $W$ is not an $L$-space.\newline
(3) $W$ admits a co-orientable taut foliation.
\end{conj}

Combining Theorem \ref{new thm} with \cite[Theorem 4.7]{R} as well as the fact that an $L$-space cannot have a co-orientable taut foliation \cite{OS}, we can update \cite[Corollary 1.3]{CD} to the following

\begin{cor}Let $M$ be the exterior of a knot in an integral
homology $3$-sphere such that $M$ fibers over the circle with
pseudo-Anosov monodromy. If the  Alexander polynomial of $M$
has a simple root on the unit
circle, then there is a real number $a>0$ such that for every rational slope $r\in (-a,a)$ the Dehn filling $M(r)$ satisfies Conjecture \ref{conj}.
\end{cor}

Similarly combining Theorem \ref{new thm} with \cite[Theorem 1.1]{LR}
we have

\begin{cor}Let $M$ be the exterior of a nontrivial knot in $S^3$.
 If the  Alexander polynomial of $M$
has a simple root on the unit
circle, then there is an $a>0$ such that for every rational slope $r\in (-a,a)$ the Dehn filling $M(r)$ satisfies Conjecture \ref{conj}.
\end{cor}

Now we proceed to  prove  Theorem \ref{new thm}.
From now on we assume that $M$ is the exterior of a knot in
an integral homology $3$-sphere such that $M$ is irreducible.

The  main new input is a quick application of some results from \cite{H} which we recall  now.
Let $F$ be a Seifert surface of genus $g$ for $M$ and let $F\times [-1,1]$ be a product  neighborhood of $F$ in $M$ so that $F=F\times \{0\}$. If $\{x_i\}_{1\leq i\leq 2g}$ is a basis for $H_1(F; \z)$, let $x_i^+$  denote the pushoff of $x_i$ in $F\times\{1\}$. Define the linking matrix $V$ by $V_{ij} =lk(x_i,x_j^+)$.
The symmetrized Alexander matrix for $M$  is the matrix
$$A(t) = t^{1/2}V -t^{-1/2}V^T$$
and $\D(t)={\rm det}A(t)$ is the Alexander polynomial of $M$.
Let $B(t) =(t^{-1/2}-t^{1/2})A(t)$. The
complex values $t\ne \pm 1$ for which $B(t )$ is singular are exactly the roots of the Alexander polynomial $\D(t)$.

Identify  $SU(2)$ with the set of unit quaternions and identify $U(1)$ with
 the unit circle in the complex plane. If $t\in U(1)$, then $B(t)$ is a Hermitian matrix and hence has only real eigenvalues.
 The equivariant knot signature of $M$, denoted by ${\rm Sign}B(t^2)$,
is the function from $U(1)$ to $\z$ taking $t$ to the number of positive eigenvalues minus the number of negative eigenvalues for $B(t^2)$, counted with multiplicity. This function is independent of the choice of
$F$, $\{x_i\}$, and the product neighborhood of $F$. The relationship between $B(t)$ and the Alexander matrix $A(t)$ implies
that ${\rm Sign}B(t^2)$ is continuous in $t\in U(1)$ except possibly at square roots of roots of the Alexander polynomial.

For each $0 <\a <\pi$, let $\r_\a:\pi_1(M)\ra SU (2)$ be
 the abelian representation determined by
  $\r(\m)=e^{i\a}$. The following results were obtained in \cite{H}.

  \begin{thm} \label{thm from H}
   {\rm (1) (\cite[Theorem 1]{H})} If  $e^{i2\a}$ is a root of $\D(t)$
       for which  the right and left hand limits $$\lim_{\b\ra \a^{\pm}}{\rm Sign}B(e^{i2\b})$$ do not agree,  then there is a continuous family of irreducible
$SU (2)$  representations of $\pi_1(M)$ limiting to $\r_\a$.
\newline
{\rm (2) (\cite[Corollary 2]{H})}
If $e^{i2\a}$ is a root of $\D(t)$ of odd multiplicity, then the condition in part (1) holds and thus there is a continuous
family of irreducible $SU (2)$ representations of $\pi_1(M )$ limiting to $\r_\a$.
\newline
{\rm (3) (\cite[Corollary 3]{H})} Suppose that
 $e^{i2\a}$ is a root of $\D(t)$ such that as $t\in U(1)$  moves through the value $e^{i\a}$, all eigenvalues
of $B(t^2)$  touching zero cross zero transversely, and all do so in the same direction.
Then all of the $SU(2)$ irreducible representations  $\{\r_s\}$ near $\r_\a$ {\rm (}provided by part {\rm (1))} send $\l$ to $e^{i\s(\r_s)}$ for some small $\s(\r_s)\ne 0$.
\end{thm}

What we need in this paper is the following special consequence of Theorem \ref{thm from H}.

\begin{cor}\label{cor needed}
If  $e^{i2\a}$ is a simple root of $\D(t)$, then there is a continuous
family of irreducible $SU (2)$ representations $\{\r_s\}$ of $\pi_1(M )$ limiting to $\r_\a$. Moreover all of these  $\r_s$  near $\r_\a$  send $\l$ to $e^{i\s(\r_s)}$ for some small $\s(\r_s)\ne 0$,
\end{cor}

\pf The first assertion is immediate by part (2) of Theorem \ref{thm from H}. To get  the second assertion, let $t_0=e^{i\a}$ and we have
$${\rm det} B(t^2)=(t^{-1}-t)^{2g}\D(t^2)=(t-t_0)f(t)$$
where $f(t)$ is a  holomorphic function
such that $f(t_0)\ne 0$.
The product rule for derivative shows that the derivative of ${\rm det} B(t^2)$ at $t_0$ is not zero.
As ${\rm det} B(t^2)$ is a product of its eigenvalues
$\l_1(t),\cdots,\l_{2g}(t)$ for which we may assume that $\l_1(t_0)=0$
and $\l_j(t_0)\ne 0$ for all $1<j\leq 2g$, applying  the  product rule
for derivative again we see that
$\l_1(t)$ has nonzero derivative at $t_0$.
So $\l_1(t)$ cross zero transversely as $t\in U(1)$  moves through the value $e^{i\a}$.
The second assertion now follows from part (3) of Theorem \ref{thm from H}.
\qed

Now let $R(M)=Hom(\pi_1M, SL_2(\c)$ be the  $SL_2(\c)$ representation variety  of $M$ and $X(M)$ the corresponding character variety.
Recall that the character  $\chi_\r\in X(M)$ of a representation $\r\in R(M)$ is the function $\chi_\r:\pi_1M\ra \c$ defined by
$\chi_\r(\g)=trace(\r(\g))$ for $\g\in \pi_1M$.
Two irreducible representations in $R(M)$ are conjugate
if and only if they have the same character.
From now on we consider $SU(2)$ as a subgroup of $SL_2(\c)$.
So the abelian representation $\r_\a$ sends $\m$ to {\small $\left(\begin{array}{cc}
e^{\i\a}&0\\0&e^{-i\a}\end{array}\right)$}.
The following results were shown in  \cite{HPP}.

\begin{thm}\label{thm,HPP}
Suppose $e^{i2\a}$ is a simple root of the Alexander polynomial $\D(t)$ of $M$. \newline
{\rm (1) (\cite[Theorem 1.2]{HPP}) }The character $\chi_{\r_\a}$ of the abelian representation $\r_\a$ is contained in a unique algebraic component
$X_0$ of $X(M)$ which contains characters of irreducible representations
and is a smooth point of $X_0$.
\newline
{\rm (2) (\cite[Theorem 1.1]{HPP}) }The complex dimension of $X_0$ is one.
\newline
{\rm (3) ( \cite[Corollary 1.4]{HPP})}
The character  $\chi_{\r_\a}$ is a middle point
of a smooth arc of real valued characters $\{\chi_{t}; t\in (-\e,\e)\}$
in $X_0$ such that $\chi_0=\chi_{\r_\a}$, $\chi_t$ is the character of an irreducible representation for  $t\ne 0$. Moreover $\chi_t$ is the character of a representation into $SU(2)$ for  $t> 0$ and $SU(1,1)$
 for  $t < 0$.
\end{thm}

Now consider a continuous family of irreducible  $SU(2)$
 representations $\{\r_s\}$  limiting to the abelian representation $\r_\a$, provided by Corollary \ref{cor needed}.
Since the character of an irreducible representation cannot be equal to
the character  of a reducible representation,
$\{\chi_{\r_s}\}$ is a nonconstant continuous family of characters
 limiting to  $\chi_{\r_\a}$.
Hence by part (1) of Theorem \ref{thm,HPP}, $\chi_{\r_s}\in  X_0$
for all $\chi_{\r_s}$ sufficiently close to $\chi_{\r_\a}$.

Let $f_\l:X_0\ra \c$ be the function defined by
$f_\l(\chi_\r)=\chi_\r(\l)-2$. Then $f_\l$ is a regular
function on the irreducible variety $X_0$.
Note that $\r_\a(\l)=I$ and by Corollary \ref{cor needed} $\r_s(\l)\ne I$ for all $\r_s$
sufficiently close to $\r_\a$, where $I$ is the identity matrix of $SL_2(\c)$. Since $\r_s(\l)\in SU(2)$, its trace cannot be $2$.
It follows that the function $f_\l$ is non-constant on $X_0$.
Since $dim_\c X_0=1$ by part (1) of Theorem \ref{thm,HPP},
any regular function  on $X_0$ is either a constant function or
has  finitely many zero points.
Hence the function $f_\l$ can have only finitely many
zeros in $X_0$. In particular
we may assume that $f_\l$ is never zero valued on the curve
$\{\chi_t; t\in (-\e,0)\}$, provided by part (3) of Theorem \ref{thm,HPP}
(by choosing smaller $\e>0$ if necessary).
The proof of \cite[Corollary 1.4]{HPP} given in \cite[Section 5]{HPP}
actually shows that there is a smooth path of irreducible
$SU(1,1)$ representations $\{\r_t; t\in (-\e,0)\}$ of $\pi_1(M)$ limiting to $\r_\a$ as $t\ra 0$ such that $\chi_t=\chi_{\r_t}$.
So for  $\r_t$, $t<0$, we have  $\r_t(\l)\ne I$.

Recall that $SU(1,1)$ is the subgroup of $SL_2(\c)$ consisting
 of matrices  {\small $\left\{  \left(\begin{array}{cc}
x&y\\\bar y&\bar x\end{array}\right);~~ x\bar x-y\bar y=1 \right\}$},
which is conjugate to $SL_2(\mathbb R)$ by the element
 {\small $\left(\begin{array}{cc}
1/\sqrt2&i/\sqrt2\\i/\sqrt2&1/\sqrt2\end{array}\right)$}.
Hence the path of irreducible $SU(1,1)$  representations
\{$\r_t; t\in (-\e,0)$\} given above is conjugate to a path of irreducible
$SL_2(\mathbb R)$ representations $\{\r_t'; t\in (-\e,0)\}$
limiting to the abelian representation $\r'_\a$ as $t\ra 0$,
where $\r_\a'$ sends $\m$ to {\small $\left(\begin{array}{cc}
\cos\a&\sin\a\\-\sin\a&\cos\a\end{array}\right)$}.
Of course $\r_t'(\l)\ne I$ for $t\in (-\e,0)$.

Now the argument in \cite[Section 7]{CD} shows that the existence
of a  path of $SL_2(\mathbb R)$ representations $\r_t'$ as given
in the preceding paragraph will imply the conclusion
of Theorem \ref{new thm}. More concretely, this path $\r_t'$ will
lift to a path of $\widetilde{SL_2(\mathbb R)}$ representations
$\widetilde{\r_t'}$ of $\pi_1(M)$, where $\widetilde{SL_2(\mathbb R)}$
 is the universal covering group of $SL_2(\mathbb R)$, and moreover
there is an $a>0$  such that for each slope $r\in (-a,0)\cup (0,a)$
some $\widetilde{\r_t'}$ will factor through $\pi_1(M(r))$
yielding a nontrivial representation $\pi_1(M(r))\ra \widetilde{SL_2(\mathbb R)}$.
Since $M$ can have at most three Dehn fillings yielding reducible manifolds by \cite{GL}, we may assume that $a>0$ has been chosen
so that $M(r)$ is irreducible for each  slope  $r\in (-a,0)\cup (0,a)$.
Hence $\pi_1(M(r))$ is left-orderable,  by \cite[Theorem 1.1]{BRW}.

\def\bysame{$\underline{\hskip.5truein}$}

\end{document}